\documentclass[12pt,draft]{amsart}
\usepackage{amsmath,amsthm,latexsym,amscd,amsbsy,amssymb}

\chardef\bslash=`\\ 

\makeatletter
\def\verbatim{\interlinepenalty\@M \@verbatim
  \leftskip\@totalleftmargin\advance\leftskip2pc
  \frenchspacing\@vobeyspaces \@xverbatim}
\makeatother
\newtheorem{thm}{Theorem}[section]
\newtheorem{cor}[thm]{Corollary}

\newtheorem{pro}[thm]{Proposition}

\newtheorem*{A}{Problem A (\cite[Problem 5]{dydakII})}
\newtheorem*{B}{Problem B (\cite[Problems 6, 7]{dydakII})}

\newtheorem*{D}{Problem A$^{\prime}$}
\newtheorem*{E}{Problem B$^{\prime}$}

\theoremstyle{definition}

\theoremstyle{remark}

\numberwithin{equation}{section}

\font\f=msbm10

\begin{document}


\title[Compactifications and universal
spaces in extension theory]
{Compactifications and universal spaces
in extension theory}
\author{Alex Chigogidze}
\address{Department of Mathematics and Statistics,
University of Saskatche\-wan,
McLean Hall, 106 Wiggins Road, Saskatoon, SK, S7N 5E6,
Canada}
\email{chigogid@math.usask.ca}
\thanks{Author was partially supported by NSERC
research grant.}
\keywords{compactification, universal space,
cohomological dimension}
\subjclass{Primary: 55M10; Secondary: 54F45}

\begin{abstract}We show that for each countable simplicial complex
$P$ the following conditions are equivalent:
\begin{itemize}
\item
$P \in AE(X)$ iff $P \in AE(\beta X)$ for any space $X$.
\item 
There exists a $P$-invertible map of a metrizable
compactum $X$
with\\
$P \in AE(X)$ onto the Hilbert cube.
\end{itemize}
\end{abstract}

\maketitle
\markboth{A.~Chigogidze}{Compactifications and universal spaces
in extension theory}


\section{Introduction}\label{S:intro}
The following two problems are central in extension theory
\cite{drady} (see also \cite{chihand}).

\begin{A}\label{P:comp}
Characterize CW-complexes $P$ such that for any space $X$ with
$P \in AE(X)$ there exist a compactification $bX$ of $X$ such
that $P \in AE(bX)$.
\end{A}

\begin{B}\label{P:universal}
Characterize CW-complexes $P$ such that the class 
\[\{ X \colon X \;\text{is a separable metrizable space with}\;
P \in AE(X) \}\]
\noindent has a universal space (compactum).
\end{B}

The first example of a space with $K(\text{\f Z},n) \in AE(X)$ and
$K(\text{\f Z},n) \notin AE(\beta X)$ (in other words
$\dim_{\text{\f Z}}X < \dim_{\text{\f Z}}\beta X$) was
constructed by A.~N.~Dranishnikov \cite{dracomp}
(see also \cite{kar}). There even exist spaces with
$K(\text{\f Z},n) \in AE(X)$ and $K(\text{\f Z},n)
\notin AE(bX)$ for any compactification $bX$ of $X$
\cite{dydakwalshcomp}. It follows from
\cite[Corollary 2.5(b) and its proof]{dydakII} that for a
finitely dominated complex $P$
the class indicated in problem B has a universal compactum.
It is important to emphasize that a universal compactum in
this case can be produced (see the proof
of \cite[Theorem 2.3]{dydakII}) as the
domain of a $P$-invertible map the
range of which is the Hilbert cube $\text{\f I}^{\omega}$.
Various results suggest and several authors have observed
(see, for instance, \cite[Remark]{dydakmogil},
\cite[p.1657]{dydakII}) that there seems to be a connection
between the existence of ``dimension" preserving compactifications
and the existence of universal elements in the class of metrizable
compacta of a given ``dimension". 

Below we consider stronger versions of the above problems. 

\begin{D}\label{P:stone}
Characterize connected locally compact simplicial complexes
\footnote{We prefer to work with connected locally compact simplicial
complexes because they are Polish ANR's and consequently 
results of \cite{book}, \cite{chicoho} apply.} $P$ such that
$P \in AE(X)$
iff $P \in AE(\beta X)$ for any space $X$.
\end{D}

\begin{E}\label{P:invert}
Characterize connected locally compact simplicial
complexes $P$ such that there exists a $P$-invertible
map $f \colon X \to \text{\f I}^{\omega}$ where $X$ is
a metrizable compactum with $P \in AE(X)$.
\end{E}

Below we show (Theorem \ref{T:main})
that problems A$^{\prime}$ and B$^{\prime}$ are equivalent.


\section{Results}\label{S:results}
All spaces are assumed to be Tychonov (i.e. completely
regular and Hausdorff). All maps are continuous.
$\text{\f I}$ denotes a closed interval. All simplicial
complexes are connected and locally compact. If $X$ is a
normal space we say that $P$ is an absolute extensor of $X$
(and write $P \in AE(X)$) if for each closed subspace $Y$
of $X$ any map $f \colon Y \to P$ has a continuous extension
$\bar{f} \colon X \to P$. An extension of this concept
for non-normal spaces has been given in
\cite[Definition 3.1]{chicoho}. A map $f \colon X \to Y$ is
$P$-invertible if for any map $g \colon Z \to Y$
with $P \in AE(Z)$
there exists a map $h \colon Z \to X$ such that $g = fh$.

\begin{thm}\label{T:main}
Let $P$ be a Polish $ANR$-space. Then the following
statements are equivalent:
\begin{itemize}
\item[(a)]
$P \in AE(\beta X)$ whenever $X$ is a space with
$P \in AE(X)$.
\item[(b)]
$P \in AE(\beta X)$ whenever $X$ is a normal space with
$P \in AE(X)$. 
\item[(c)]
$P \in AE\left(\beta\left(\oplus \{ X_{t} \colon t \in T\}
\right)\right)$ whenever $T$ is an arbitrary indexing set
and $X_{t}$, $t \in T$, is a separable metrizable space
with $P \in AE(X_{t})$. 
\item[(d)]
$P \in AE\left(\beta\left(\oplus \{ X_{t} \colon t \in T\}
\right)\right)$ whenever $T$ is an arbitrary indexing set
and $X_{t}$, $t \in T$, is a Polish space with $P \in AE(X_{t})$. 
\item[(e)]
There exists a $P$-invertible map $f_{P} \colon K_{P} \to
\text{\f I}^{\omega}$ where $K_{P}$ is a metrizable compactum
with $P \in AE(K_{P})$.
\end{itemize}
\end{thm} 
\begin{proof}
Implications (a) $\Longrightarrow$ (b) $\Longrightarrow$ (c)
$\Longrightarrow$ (d) are trivial. Proof of implication (d)
$\Longrightarrow$ (e) follows the proof of \cite[Proposition
5.3]{chicoho}. Let $\mathcal A$ denote the set of all maps
$\{ r_{t} \colon t \in T\}$ such that domain
$\operatorname{Dom}(r_{t})$ is a Polish subspace of
$\text{\f I}^{\omega}$,
$P \in AE\left( \operatorname{Dom}(r_{t})\right)$ and
$\operatorname{Ran}(r_{t}) \subseteq \text{\f I}^{\omega}$.
Let $Y = \oplus\{ \operatorname{Dom}(r_{t}) \colon t \in T\}$.
Clearly, $P \in AE(Y)$. Consider also the map $r \colon Y
\to \text{\f I}^{\omega}$ which coincides with $r_{t}$ on
$\operatorname{Dom}(r_{t})$ for each $t \in T$. Let
$\bar{r} \colon \beta Y \to \text{\f I}^{\omega}$ be the unique
continuous extension of $r$ to the Stone-\v{C}ech compactification
$\beta Y$ of $Y$. By (d), $P \in AE(\beta Y)$. By
\cite[Theorem 4.4]{chicoho} and by the {\em compactness} of
$\beta Y$, the latter is the limit space of a Polish
spectrum ${\mathcal S} = \{ Y_{\alpha}, q_{\alpha}^{\beta},
A\}$ consisting of metrizable {\em compacta} $Y_{\alpha}$
(compactness of $Y_{\alpha}$
follows from the fact that $q_{\alpha}(\beta Y)$ is dense
in $Y_{\alpha}$, according to assmption made in
\cite[p. 201]{chicoho})
with $P \in AE(Y_{\alpha})$. Write
$\text{\f I}^{\omega} =
\prod\{ \text{\f I}_{n} \colon n \in \omega\}$ where
$\text{\f I}_{n}$, $n \in \omega$, denotes a copy
of $\text{\f I}$.
Let also $\pi_{n} \colon \text{\f I}^{\omega}
\to \text{\f I}_{n}$, $n \in \omega$, denote the
corresponding projection. Since the spectrum $\mathcal S$
is factorizing, for each $n \in \omega$ there
exist an index $\alpha_{n} \in A$ and a map
$s_{n} \colon Y_{\alpha_{n}}
\to \text{\f I}_{n}$ such that
$\pi_{n}\bar{r} = s_{n}q_{\alpha_{n}}$,
where $q_{\alpha_{n}} \colon \beta Y \to Y_{\alpha_{n}}$ is the
$\alpha_{n}$-th limit projection of the spectrum $\mathcal S$.
Since $\mathcal S$ is a Polish spectrum
(see \cite[page 201]{chicoho}) there exists an
index $\alpha \in A$ such that $\alpha \geq \alpha_{n}$
for each $n \in \omega$. Next consider the map
\[ s = \triangle\{ s_{n}q_{\alpha_{n}}^{\alpha} \colon n \in \omega \}
\colon Y_{\alpha} \to \prod\{ \text{\f I}_{n}
\colon n \in \omega\} ,\]
where
$q_{\alpha_{n}}^{\alpha} \colon Y_{\alpha} \to Y_{\alpha_{n}}$,
$n \in \omega$, denotes the corresponding projection
of the spectrum $\mathcal S$. It is easy to see that
$\bar{r} = sq_{\alpha}$, where
$q_{\alpha} \colon \beta Y \to Y_{\alpha}$ is the
$\alpha$-th limit projection of the spectrum $\mathcal S$.
It now suffices to let $K_{P} = Y_{\alpha}$ and $f_{P} = s$.
Let us show that $s \colon Y_{\alpha} \to \text{\f I}^{\omega}$
is indeed $P$-invertible. Since the spaces
$Y_{\alpha}$ and $\text{\f I}^{\omega}$ are Polish
(even compact and metrizable), it suffices
(according to \cite[Proposition 5.2]{chicoho})
to consider only Polish spaces $Z$ in the definition of
$P$-invertibility given above. Indeed, let
$g \colon Z \to \text{\f I}^{\omega}$ be a map
defined on a Polish space
$Z$ with $P \in AE(Z)$. We may as well assume that
$Z \subseteq \text{\f I}^{\omega}$. By the definition of
$\mathcal A$,
there is an index $t \in T$ such that $r_{t} = g$.
Let $i_{t} \colon \operatorname{Dom}(r_{t}) \to Y$ denote the
corresponding embedding. Clearly, $r_{t} = \bar{r}i_{t}$.
Then the composition
$h = q_{\alpha}i_{t} \colon Z \to Y_{\alpha}$ lifts
the map $g$, i.e. $sh = g$.

(e) $\Longrightarrow$ (a). As in the proof of \cite[Theorem 5.13]
{chicoho} (see also \cite[Section 6.2]{book} where the case
$P = S^{n}$ is considered) one shows that for any uncountable
cardinal number $\tau$ there exists a $P$-invertible map
$f = f_{P,\tau} \colon K_{P,\tau} \to \text{\f I}^{\tau}$, where
$K_{P,\tau}$ is a compactum of weight $\tau$ such that $P \in
AE\left( K_{P,\tau}\right)$ and $\text{\f I}^{\tau}$ denotes
the Tychonov cube of weight $\tau$ (for $\tau = \omega$ the
existence of such a map is guaranteed by condition (e);
note also that an $AE(0)$-space of countable weight
is Polish \cite[Corollary 6.4.5]{book}). Consider
now a space $X$ with $P \in AE(X)$ and choose $\tau$ large
enough so that $\beta X$ can be identified with a subspace
of $\displaystyle \text{\f I}^{\tau}$. Since the map $f$ is
$P$-invertible there exists a map $g \colon X \to K_{P,\tau}$
such that $fg = \operatorname{id}_{X}$. Since $K_{P,\tau}$ is
compact, the map $g$ admits a continuous extension $\bar{g} \colon \beta X
\to K_{P,\tau}$. Since $fg = \operatorname{id}_{X}$ and since
$\bar{g}|X = g$ it follows that $f\bar{g} =
\operatorname{id}_{\beta X}$. In this situation it can easily
be seen that $\bar{g}$ is an embedding. In other words,
$\bar{g}(\beta X)$ is a topological copy of $\beta X$. Finally,
since $P \in AE(K_{P,\tau})$ it follows that $P \in AE(\beta X)$.
\end{proof}

\begin{cor}\label{C:dominated}
Let $P$ be a finitely dominated connected locally compact complex. 
Then the following conditions are equivalent for any space $X$:
\begin{enumerate}
\item
$P \in AE(X)$.
\item
$P \in AE(\beta X)$.
\end{enumerate}
\end{cor}
\begin{proof}
The implication (2) $\Longrightarrow$ (1) follows from
\cite[Proposition 6.8]{chicoho}. Let us prove the
implication (1) $\Longrightarrow$ (2). According to
Theorem \ref{T:main}
it suffices to construct a $P$-invertible map
$f_{P} \colon K_{P} \to \text{\f I}^{\omega}$ where $K_{P}$
is a metrizable compactum with $P \in AE(K_{P})$. Since $P$
is finitely dominated, there exist a finite complex $L$ and
two maps $u \colon P \to L$ and $d \colon L \to P$ such
that $du \simeq \operatorname{id}_{P}$. By 
\cite[Theorem 2.3]{dydakII} there exists a compactum
$K_{P}$ and a map $f_{P} \colon K_{P} \to \text{\f I}^{\omega}$
with the following properties:
\begin{itemize}
\item[(a)]
For each map $g \colon Z \to \text{\f I}^{\omega}$,
defined on a separable
metrizable space with $P \in AE(Z)$, there exists a
map $h \colon Z \to K_{P}$ such that $f_{P}h = g$.
\item[(b)]
For each map $\varphi \colon C \to P$, where $C$ is a closed
subset of $K_{P}$, there exists a map
$\varphi^{\prime} \colon K_{P} \to L$
such that $\varphi^{\prime}|C \simeq u\varphi$.
\end{itemize}
Observe that, by (b), $d\varphi^{\prime}|C \simeq du\varphi \simeq \varphi$.
Consequently, by the Homotopy Extension Theorem,
the map $\varphi$ has a continuous extension over $K_{P}$.
This in turn means that $P \in AE(K_{P})$. By (a) and 
\cite[Proposition 5.2]{chicoho}, the map $f_{P}$ is $P$-invertible.
\end{proof}

\begin{pro}\label{P:finitetype}
Let $P$ be a connected locally compact simplicial
complex of finite type with a finite fundamental group.
Then the following conditions are equivalent for any space
$X$ and any integer $n \geq 2$:
\begin{itemize}
\item[(a)]
$P \vee S^{n} \in AE(X)$.
\item[(b)]
$P \vee S^{n} \in AE(\beta X)$.
\end{itemize}
\end{pro}
\begin{proof}
The implication (b) $\Longrightarrow$ (a) follows from
\cite[Proposition 6.8]{chicoho}.

In order to prove the implication (a) $\Longrightarrow$ (b)
it suffices to show that $P \vee S^{n} \in AE(\beta X)$
for each normal space $X$ with $P \vee S^{n} \in AE(X)$
(see Theorem \ref{T:main}(b)). Let $\varphi \colon F \to P
\vee S^{n}$ be a map defined on a closed subset $F$ of
$\beta X$. Since $P \vee S^{n}$ is an $ANR$, there exists
an extension $\psi \colon \operatorname{cl}_{\beta X}V \to
P \vee S^{n}$ of $\varphi$, where $V$ is an open neighborhood
of $F$ in $\beta X$. Clearly $V \cap X \neq \emptyset$ and
$\operatorname{cl}_{\beta X}(\left(\operatorname{cl}_{X}
(V \cap X)\right) = \operatorname{cl}_{\beta X}V$. Since
$P \vee S^{n} \in AE(X)$ there exists a map $f \colon X \to
P \vee S^{n}$ such that $f|\operatorname{cl}_{X}(V \cap X) =
\psi|\operatorname{cl}_{X}(V \cap X)$. An argument similar to
\cite[Proof of Lemma 4.1]{caldersiegel} shows that
$f$ is homotopic to a map $g \colon X \to P \vee S^{n}$
such that $\operatorname{cl}\left( g(X)\right)$ is compact.
Consequently $g$ has a continuous extension $\bar{g} \colon \beta X
\to P \vee S^{n}$ onto the whole $\beta X$. Now consider the 
two maps $\psi$ and $\bar{g}|\operatorname{cl}_{\beta X}V$.
Their restrictions $\psi |\operatorname{cl}_{X}(V\cap X)$ and $\bar{g}|\operatorname{cl}_{X}(V\cap X)$ are homotopic.
Since $X$ is normal it follows that
$\operatorname{cl}_{\beta X}V = \operatorname{cl}_{\beta X}
(\left(\operatorname{cl}_{X}(V \cap X)\right) = \beta\left(
\operatorname{cl}_{X}(V \cap X)\right)$. By
\cite[Theorem 4.2]{caldersiegel}, the restriction operator
provides a bijection of homotopy classes $\left[\beta \left(
\operatorname{cl}_{X}(V \cap X)\right) ,P \vee S^{n}\right]$
and $\left[ \operatorname{cl}_{X}(V \cap X),
P \vee S^{n}\right]$. Consequently, $\psi \simeq
\bar{g}|\operatorname{cl}_{\beta X}V$. By the homotopy
extension theorem (recall that $\bar{g}|
\operatorname{cl}_{\beta X}V$ has an extension
$\bar{g}$ onto $\beta X$ and that $P \vee S^{n}$ is
an $ANR$), $\psi$ also has an extension onto $\beta X$
which serves as an extension of the originally given
map $\varphi$.
\end{proof}

Finally I would like to thank the referee for
helpful suggestions which led to a number of improvments
of the
original version of these notes.

\end{document}